# Exact results for deterministic cellular automata traffic models

Henryk Fukś

*The Fields Institute for Research*
*in Mathematical Sciences*
*Toronto, Ontario M5T 3J1, Canada*
*and*
*Department of Mathematics and*
*Statistics,*
*University of Guelph,*
*Guelph, Ontario N1G 2W1, Canada*

Email: hfuks@fields.utoronto.ca

**Abstract**

We present a rigorous derivation of the flow at arbitrary time in a deterministic cellular automaton model of traffic flow. The derivation employs regularities in preimages of blocks of zeros, reducing the problem of preimage enumeration to a well known lattice path counting problem. Assuming infinite lattice size and random initial configuration, the flow can be expressed in terms of generalized hypergeometric function. We show that the steady state limit agrees with previously published results.

## 1. Introduction

Since the introduction of the Nagel-Schreckenberg (N-S) model in 1992 [8], cellular automata became a well established method of traffic flow modeling. Comparatively low computational cost of cellular automata models made it possible to conduct large-scale real-time simulations of urban traffic in the city of Duisburg [2] and Dallas/Forth Worth [10]. Several simplified models have been proposed, including models based on deterministic cellular automata. For example, Nagel and Herrmann [9] considered deterministic version of the N-S model, while Fukui and Ishibashi [5] introduced another model (to be referred to as F-I model), which can be understood as a generalization of cellular automaton rule 184. Rule 184, one of the elementary CA rules investigated by Wolfram [13], had been later studied in detail as a simple model of surface growth [6], as well as in the context of density classification problem [3]. It is one of the only two (symmetric) non-trivial elementary rules conserving the number of active sites [1], and, therefore, can be interpreted as a rule governing dynamics of particles (cars). Particles (cars) move to the left if their right neighbor site is empty, and do not move if the right neighbor site is occupied, all of them moving simultaneously at each discrete time step. Using terminology of lattice stochastic processes, rule 184 can be viewed as a discrete-time version of totally asymmetric simple exclusion process. Further generalization of the F-I model has been proposed in [4].

In all traffic models, the main quantity of interest is the average velocity of cars, or the average flow, defined as a product of the average velocity and the density of cars. The graph of the flow as a function of density is called a fundamental diagram, and is typically studied in the steady state ($t \to \infty$). For the F-I model, steady-state fundamental diagram can be obtained using mean-field argument [5], as well as by statistical mechanical approach [11] or by studying the time evolution of inter-car spacing [12]. In general, little is known about non-equilibrium properties of the flow. In [3], we investigated dynamics of rule 184 and derived expression for the flow at arbitrary time, assuming that the initial configuration (at $t = 0$) was random, using the concept of defects and analyzing the dynamics of their collisions. In what follows, we shall generalize results of [3] for the deterministic F-I traffic flow model and derive the expression for the flow at arbitrary time. The derivation employs regularities of preimages of blocks of zeros, reducing the problem of preimage enumeration to a well known combinatorial problem of lattice path counting. Assuming infinite lattice size and random initial configuration, the flow can then be expressed in terms of generalized hypergeometric function. We will, unlike in [3], explore regularities of preimages using purely algebraic methods, i.e., without resorting to properties of spatiotemporal diagrams and dynamics of defects.

## 2. Deterministic traffic rules

Deterministic version of the F-I traffic model is defined on one-dimensional lattice of $L$ sites with periodic boundary conditions. Each site is either occupied by a vehicle, or empty. The velocity of each vehicle is an integer between 0 and $m$. If $x(i, t)$ denotes the position of the $i$th car at time $t$, the position of the next car ahead at time $t$ is $x(i + 1, t)$. With this notation, the system evolves according to a synchronous rule given by

$$x(i, t + 1) = x(i, t) + v(i, t), \tag{1}$$

where

$$v(i, t) = \min\left(x(i + 1, t) - x(i, t) - 1, m\right) \tag{2}$$

is the velocity of car $i$ at time $t$. Since $g = x(i + 1, t) - x(i, t) - 1$ is the gap (number of empty sites) between cars $i$ and $i + 1$ at time $t$, one could say that each time step, each car advances by $g$ sites to the right if $g \leq m$, and by $m$ sites if $g > m$. When $m = 1$, this model is equivalent to elementary cellular automaton rule 184, for which a number of exact results is known [6, 3].

The main quantities of interest in this paper will be the average velocity of cars at time $t$ defined as

$$\overline{v}(t) = \frac{1}{N} \sum_{i=1}^{N} v(i, t), \tag{3}$$

and the average flow $\phi(t) = \rho \overline{v}(t)$, where $\rho = N/L$ is the density of cars. In what follows, we will assume that at $t = 0$ the cars are randomly distributed on the lattice. When $N \to \infty$, this corresponds to a situation when sites are occupied by a car with probability $\rho$, or are empty with probability $1 - \rho$.

In general, if $N_k(t)$ is the number of cars with velocity $k$, we have

$$\overline{v}(t) = \frac{1}{N} \sum_{k=1}^{m} k N_k(t). \tag{4}$$

When $k < m$, $N_k(t)$ is just the number of blocks of type $10^k 1$, where $0^k$ denotes $k$ zeros. This means that a probability of an occurrence of the block $10^k 1$ at time $t$ can be written as $P_t(10^k 1) = N_k/L$. Similarly, for $k = m$, $P_t(10^m) = N_m(t)/L$. As a consequence, equation (4) becomes

$$\overline{v}(t) = \sum_{k=1}^{m-1} \frac{k P_t(10^k 1)}{\rho} + \frac{m P_t(10^m)}{\rho} \tag{5}$$

We will now demonstrate that in the deterministic F-I model with maximum speed $m$ the average flow depends only on one block probability. More precisely, we shall prove the following:

**Proposition 1.** *In the deterministic F-I model with the maximum speed $m$, the average flow $\phi_m(t)$ is given by*

$$\phi_m(t) = 1 - \rho - P_t(0^{m+1}). \tag{6}$$

To prove this proposition by induction, we first note that for $m = 1$ equation (5) gives $\phi_1(t) = P_t(10)$. Using consistency condition for block probabilities $P_t(10) + P_t(00) = P_t(0) = 1 - \rho$, we obtain $\phi_1(t) = 1 - \rho - P(00)$, which verifies (6) in the $m = 1$ case. Now assume that (6) is true for some $m = n - 1$ (where $n > 1$), and compute $\phi_n(t)$:

$$\begin{aligned}
\phi_n(t) &= n P_t(10^n) + \sum_{j=1}^{n-1} j P(10^j 1) = \\
&= n P_t(10^n) + (n-1) P_t(10^{n-1} 1) + \sum_{j=1}^{n-2} j P(10^j 1) \\
&= (n-1)[P_t(10^{n-1} 1) + P_t(10^n)] + P_t(10^n) + \sum_{j=1}^{n-2} j P(10^j 1)
\end{aligned}$$

Using consistency condition $P_t(10^{n-1} 1) + P_t(10^n) = P_t(10^{n-1})$ we obtain

$$\phi_n(t) = P_t(10^n) + (n-1) P_t(10^{n-1}) + \sum_{j=1}^{n-2} j P(10^j 1) = P_t(10^n) + \phi_{m-1}(t)$$

Taking into account that $P_t(10^n) = P_t(0^n) - P_t(0^{n+1})$ (which, again, is just a consistency condition for block probabilities), and using (6) to express $\phi_{m-1}(t)$, we finally obtain

$$\phi_m(t) = 1 - \rho - P_t(0^{m+1}). \tag{7}$$

This means that validity of (6) for $m = n$ follows from its validity for $m = n - 1$, concluding our proof by induction.

## 3. Enumeration of preimages of $0^{m+1}$

Proposition 1 reduces the problem of computing $\phi_m(t)$ to the problem of finding the probability of a block of $m+1$ zeros. In order to find this probability, we will now use the fact that the deterministic F-I model is equivalent to a cellular automaton defined as follows. Let $s(i,t)$ denotes the state of a lattice site $i$ at time $t$ (note that $i$ now labels consecutive lattice sites, *not* consecutive cars), where $s(i,t) = 1$ for a site occupied by a car and $s(i,t) = 0$ otherwise. We can immediately realize that if a site $i$ is empty at time $t$, then at time $t+1$ it can become occupied by a car arriving from the left, but not from a site further than $i - m$. Similarly, if a site $i$ is occupied, it will become empty at the next time step only and only if site $i + 1$ is empty. Thus, in general, $s(i, t+1)$ depends on $s(i-m,t), s(i-m+1,t), \ldots, s(i+1,t)$, i.e., on the state of $m$ sites to the left, one site to the right, and itself, but not on any other site, what can be expressed as

$$s(i, t+1) = f_m\Big(s(i-m,t), s(i-m+1,t), \ldots, s(i+1,t)\Big), \qquad (8)$$

where $f_m$ is called a local function of the cellular automaton. For the F-I CA, one can write explicit formula[1] for $f_m$, such as

$$f_m\Big(s(i-m,t), s(i-m+1,t), \ldots, s(i+1,t)\Big) = s(i,t) - \min\{s(i,t), 1-s(i+1,t)\}$$
$$+ \min\Big\{\max\{s(i-m,t), s(i-m+1,t), \ldots, s(i-1,t)\}, 1-s(i,t)\Big\}, \qquad (9)$$

which, using terminology of cellular automata theory, represents a rule with left radius $m$ and right radius 1. In general, after $t$ iteration of this cellular automaton rule, state of a site $s(i,t)$ depends on $s(i-mt,0), s(i-mt+1,0), \ldots, s(i+t,0)$, but not on any other sites in the initial configuration. Similarly, a block of $k$ sites $s(i,t)s(i+1,t)\ldots s(i+k)$ depends only on a block $s(i-mt,0), s(i-mt+1,0), \ldots, s(i+k+t,0)$, as schematically shown in Figure 1. We will say that $s(i-mt,0), s(i-mt+1,0), \ldots, s(i+k+t,0)$ is an $t$-step preimage of the block $s(i,t)s(i+1,t)\ldots s(i+k)$. Preimages in the F-I cellular automaton have the following property:

**Proposition 2.** *Block $a_1 a_2 a_3 \ldots a_p$ is an $n$-step preimage of a block $0^{m+1}$ if and only if $p = (n+1)(m+1)$ and, for every $k$ ($1 \leq k \leq p$)*

$$\sum_{i=1}^{k} \xi(a_i) > 0, \qquad (10)$$

*where $\xi(1) = -m$ and $\xi(0) = 1$.*

Before we present a proof of this proposition, note that it can be interpreted as follows. Let us assume that we have a block of zeros and ones of length $p$, where $p = (n+1)(m+1)$, and we want to check if this block is an $n$-step preimage of a block $0^{m+1}$. We start with a "capital" equal to zero. Now we move from the leftmost site to the right, and every time we encounter 0, we increase our capital by $m$. Every time we encounter 1, our capital decreases by 1. If we can move from $a_1$ to $a_p$ and our capital stays always larger than zero, the string $a_1 a_2 a_3 \ldots a_p$ is a preimage of $0^{m+1}$. Condition (10) can be also written as

$$\sum_{i=1}^{k} a_i < \frac{k}{m+1}, \tag{11}$$

because $\xi(x) = 1 - (m+1)x$ for $x \in \{0, 1\}$.

---
[1]Since formula (9) will not be used in subsequent calculations, we give it withot proof (which is elementary).

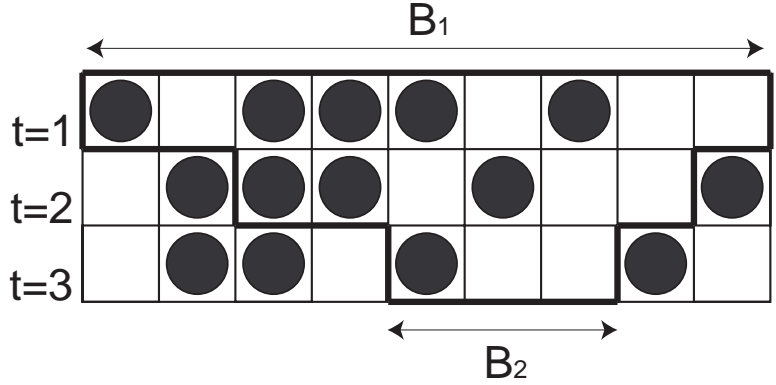

Figure 1: Fragment of a spatiotemporal diagram for the F-I rule with $m = 2$. States of nine sites during three consecutive time steps are are shown, black circles representing occupied sites. Block $B_1 = 101110100$ is a 2-step preimage of the block $B_2 = 100$. Outlined sites constitute "light cone" of the block $B_2$, meaning that the state of sites belonging to $B_2$ can depend only on sites inside the outlined region, but not on sites outside this region.

For the purpose of the proof, strings $a_1 a_2 \ldots a_p$ of length $p$ satisfying (11) for a given $m$ and for every $k \leq N$ will be called $m$-admissible strings.

**Lemma.** Let $s(1,t)s(2,t)\ldots s(p,t)$ be an $m$-admissible string. If

$$s(i, t+1) = f_m\bigl(s(i-m,t), s(i-m+1,t), \ldots, s(i+1,t)\bigr), \tag{12}$$

and if $f_m$ is a local function of the deterministic F-I model with maximum speed $m$, then $s(m+1, t+1)s(m+2, t+1)\ldots s(p-1, t+1)$ is also an $m$-admissible string.

To prove the lemma, it is helpful to employ the fact that the F-I rule conserves the number of cars. Let $0 < k < p$ and let us consider strings $S_1 = s(1,t)s(2,t)\ldots s(k,t)$ and $S_2 = s(m+1, t+1)s(2, t+1)\ldots s(k, t+1)$. If the string $s(1,t)s(2,t)\ldots s(k,t)$ is $m$-admissible, then its first $m+1$ sites must be zeros. This means that in one time step, no car can enter string $s(1,t)s(2,t)\ldots s(k,t)$ from the left. On the other hand, in a single time step, only one car (or none) can leave the string on the right hand side, i.e.,

$$\sum_{i=1}^{k} s(i,t) = \epsilon + \sum_{i=m+1}^{k} s(i, t+1) \tag{13}$$

where $\epsilon \in \{0, 1\}$. Three cases can be distinguished:

**(i)** All sites $s(k-m+1,t)s(k-m+2,t)\ldots s(k,t)$ are empty (equal to 0). Then no car leaves $S_1$, which means that $\epsilon = 0$, and

$$\sum_{i=1}^{k} s(i,t) = \sum_{i=m+1}^{k} s(i,t+1) = \sum_{i=m+1}^{k-m} s(i,t+1) < \frac{k-m}{m+1}. \tag{14}$$

The last inequality is a direct consequence of $m$-admissibility of $S_1$. Since the length of the string $S_2$ is equal to $k - m$, the above relation (which holds for arbitrary $k$) proves that $S_2$ is also $m$-admissible in the case considered.

**(ii)** Among sites $s(k-m+1,t)s(k-m+2,t)\ldots s(k,t)$ there is at least one which is occupied (equal to 1), and $s(k+1,t) = 1$. In this case, since the last site in $S_1$ is "blocked" by the car at $s(k+1,t)$, again no car can leave string $S_1$ in one time step. Therefore,

$$\sum_{i=1}^{k} s(i,t) = \sum_{i=m+1}^{k} s(i,t+1). \tag{15}$$

$m$-admissibility of $S_2$ implies

$$\frac{k+1}{m+1} > \sum_{i=1}^{k+1} s(i,t) = \sum_{i=1}^{k} s(i,t) + 1. \tag{16}$$

Combining (15) with (16) we obtain

$$\sum_{i=m+1}^{k-m} s(i,t+1) < \frac{k-m}{m+1}, \tag{17}$$

which again shows that $S_2$ is $m$-admissible.

**(iii)** Among sites $s(k-m+1,t)s(k-m+2,t)\ldots s(k,t)$ there is at least one which is occupied (equal to 1), and $s(k+1,t) = 0$. In this case, one car will leave right end of the string $S_1$, therefore

$$\sum_{i=1}^{k} s(i,t) = \sum_{i=m+1}^{k} s(i,t+1) - 1. \tag{18}$$

As before, from $m$-admissibility of $S_1$ we have

$$\sum_{i=1}^{k+1} s(i,t) = \sum_{i=1}^{k} s(i,t) < \frac{k+1}{m+1}, \tag{19}$$

hence

$$\sum_{i=m+1}^{k} s(i,t+1) = \sum_{i=m+1}^{k} s(i,t+1) - 1 < \frac{k+1}{m+1} - 1 = \frac{k-m}{m+1}, \tag{20}$$

which demonstrates that case (iii) also leads to $m$-admissibility of $S_2$, concluding the proof of our lemma.

Let us now assume that the block $B_1 = s(1,t)s(2,t)\ldots s(p,t)$ is $m$ admissible ($n$ being some fixed integer and $p = (m+1)(n+1)$). Applying the lemma to this block we conclude that $B_2 = s(m+1,t+1)s(2,t+1)\ldots s(p-1,t+1)$ is $m$-admissible as well. Applying the lemma to $B_2$ we obtain $m$-admissible block $B_3 = s(2m+1, t+2)s(2, t+2) \ldots s(p-2, t+2)$. After $n$ applications of the lemma we end up with the conclusion that the string $B_{n+1} = s(nm+1, n+1)s(nm+2, n+1) \ldots s(p-n)$ is $m$-admissible. Since the length of $B_{n+1}$ is $p - n - nm = (n+1)(m+1) - n(m+1) = m+1$, it must, to be $m$-admissible, be composed of all zeros, i.e., $B_{m+1} = 0^{m+1}$. This means that $m$-admissibility of $B_1$ is a sufficient condition for $B_1$ to be an $n$-step preimage of $0^{m+1}$. Reversing steps in the above reasoning, one can show that is is also a necessary condition.

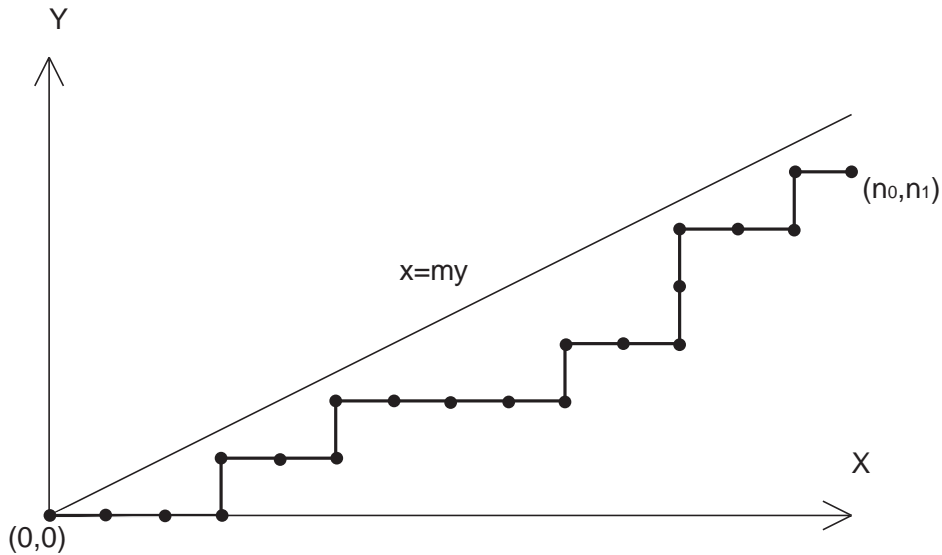

Figure 2: $m$-admissible block with $n_0$ zeros and $n_1$ ones is equivalent to a lattice path from the origin to $(n_0, n_1)$ which does not touch nor cross the line $x = my$. 0 corresponds to a horizontal segment, while 1 to a vertical segment.

## 4. Fundamental diagram

We shall now use proposition 2 to calculate $P_t(0^{m+1})$. First of all, we note that $P_t(0^{m+1})$ is equal to the probability of occurrence of $t$-step preimage of $0^{m+1}$ in the initial (random) configuration, that is

$$P_t(0^{m+1}) = \sum P_0(a), \qquad (21)$$

where the sum goes over all $t$-step preimages of $0^{m+1}$. Consider now a string which contains $n_0$ zeros and $n_1$ ones. The number of such strings can be immediately obtained if we realize that it is equal to the number of lattice paths from the origin to $(n_0, n_1)$ which do not touch nor cross the line $x = my$, as shown in Figure 2. This is a well known combinatorial problem [7], and the number of aforementioned paths equals

$$\frac{n_0 - mn_1}{n_0 + n_1} \binom{n_0 + n_1}{n_1} \qquad (22)$$

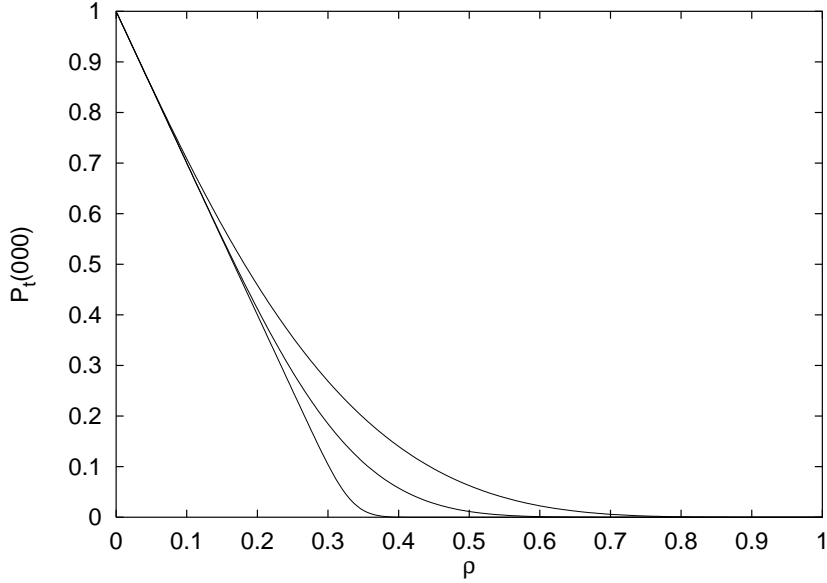

Figure 3: Graph of the probability $P_t(0^{m+1})$ as a function of $\rho$ for $m = 2$ and $t = 1$ (upper line), $t = 5$ (middle line), and $t = 100$ (lower line).

Probability of occurrence of such a block in a random configurations is, therefore,

$$\frac{n_0 - mn_1}{n_0 + n_1}\binom{n_0 + n_1}{n_1}\rho^{n_1}(1-\rho)^{n_0}, \qquad (23)$$

where $\rho = P(1)$. In a $t$-step preimage of $0^{m+1}$ the minimum number of zeros is $1 + m(t+1)$ zeros, while the maximum is $(m+1)(t+1)$ (corresponding to all zeros). Therefore, summing over all possible number of zeros $i$, we obtain

$$P_t(0^{m+1}) = \sum_{i=1+m(t+1)}^{(m+1)(t+1)} \frac{i - m[(m+1)(t+1) - i]}{(m+1)(t+1)}\binom{(m+1)(t+1)}{(m+1)(t+1) - i}$$
$$\times \rho^{(m+1)(t+1)-i}(1-\rho)^i$$

Changing summation index $j = i - m(t+1)$ we obtain

$$P_t(0^{m+1}) = \sum_{j=1}^{t+1} \frac{j}{t+1}\binom{(m+1)(t+1)}{t+1-j}\rho^{t+1-j}(1-\rho)^{m(t+1)+j}. \qquad (24)$$

Figure 3. shows a graph of $P_t(0^{m+1})$ as a function of $\rho$ for $m = 2$ and several values of $t$. We can observe that as $t$ increases, the graph becomes "sharper" at $\rho = 1/3$, eventually developing singularity (discontinuity in the first derivative) at $\rho = 1/3$. More precisely, one can show (see appendix) that

$$\lim_{t \to \infty} P_t(0^{m+1}) = \begin{cases} 1 - (m+1)\rho & \text{if } p < 1/(m+1), \\ 0 & \text{otherwise.} \end{cases} \qquad (25)$$

$P_\infty(0^{m+1})$, therefore, can be viewed as the order parameter in a phase transition with critical point at $\rho = 1/(m+1)$. Using Proposition 1 we can now find the average flow in the steady state

$$\phi_m(\infty) = \begin{cases} m\rho & \text{if } p < 1/(m+1), \\ 1-\rho & \text{otherwise,} \end{cases} \qquad (26)$$

which agrees with mean-field type calculations reported in [5] as well as with results of [11, 12].

To verify validity of the result for $t < \infty$, we performed computer simulations using a lattice of $10^5$ sites with periodic boundary conditions. The average flow has been recorded after each iteration up to $t = 100$ for three values of $\rho$: at the critical point $\rho = 1/3$ as well as below and above the critical point. The resulting plots of the flow as a function of time are presented in Figure 4. Again, the agreement with theoretical curves

$$\phi_m(t) = 1 - \rho - \sum_{j=1}^{t+1} \frac{j}{t+1} \binom{(m+1)(t+1)}{t+1-j} \rho^{t+1-j}(1-\rho)^{m(t+1)+j} \qquad (27)$$

is very good.

Without going into details, we note that the formula (27) can be also expressed in terms of generalized hypergeometric function $_2F_1$:

$$\phi_m(t) = 1 - \rho - \frac{(1-\rho)^{1+m+mt}\rho^t(1+m+t+mt)!}{(1+m+mt)(1+t)!(m+mt)!} \, _2F_1\left[\begin{array}{c} 2, -t \\ 2+m+mt \end{array}; 1 - \frac{1}{\rho}\right], \qquad (28)$$

Since fast numerical algorithms for computing $_2F_1$ exist, this form might be useful for the purpose of numerical evaluation of $\phi_m(t)$.

## 5. Conclusion

We presented derivation of the flow at arbitrary time in the deterministic F-I cellular automaton model of traffic flow. First, we showed that the flow can be expressed by the probability of occurrence of the block of $m+1$ zeros $P(0^{m+1})$. By employing regularities in preimages of blocks of zeros, we reduced the problem of preimage enumeration to the lattice path counting problem. Finally, we used the number of preimages to find $P(0^{m+1})$, which determines the flow.

We also found that the flow in the steady state, obtained by taking $t \to \infty$ limit, agrees with previously reported mean-field type calculations, meaning that in the case of the F-I model mean-field approximation gives exact results. This seems to be true not only for the F-I model, but also for many other CA rules conserving the number of active sites ("conservative" CA). For example, in [1] we reported that the third order local structure approximation, which is a generalization of simple mean-field theory incorporating short-range correlations, yields the fundamental diagram for rule 60200 (one of the 4-input "conservative" CA rules) in extremely good agreement with computer simulations. Taking this into account, we conjecture that the local structure approximation gives exact fundamental diagram for almost all "conservative" rules, excluding, perhaps, those rules for which the fundamental diagram is not sufficiently "regular" (meaning not piecewise linear). This problem is currently under investigation.

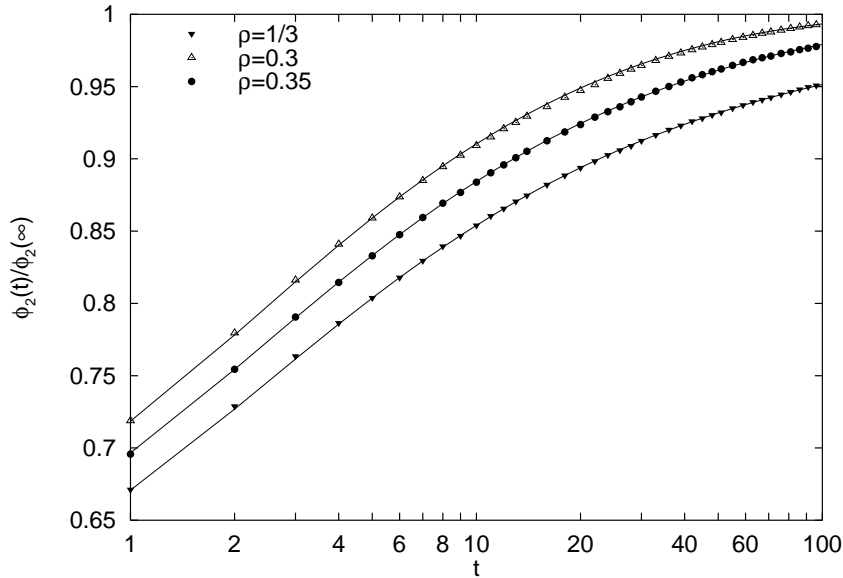

Figure 4: Plots of $\phi_2(t)/\phi_2(\infty)$ as a function of time for $\rho = 0.3$, $\rho = 1/3$, and $\rho = 0.35$ obtained from computer simulation on a lattice of $10^5$ sites. Continuous line corresponds to the theoretical result obtained using eq. (28).

## Acknowledgements

The author wishes to thank The Fields Institute for Research in Mathematical Sciences for generous hospitality and the National Sciences and Engineering Research Council of Canada for financial support.

## Appendix

In order to find the limit $\lim_{t\to\infty} P_t(0^{m+1})$ we can write eq. (24) in the form

$$P_t(0^{m+1}) = \sum_{j=1}^{t+1} \frac{j}{t+1} b(t+1-j, (m+1)(t+1), \rho), \qquad (29)$$

where

$$b(k,n,p) = \binom{n}{k} p^k (1-p)^{n-k} \qquad (30)$$

is the distribution function of the binomial distribution. Using de Moivre-Laplace limit theorem, binomial distribution for large $n$ can be approximated by the normal distribution

$$b(k,n,p) \sim \frac{1}{\sqrt{2\pi np(1-p)}} \exp \frac{-(k-np)^2}{2np(1-p)}. \qquad (31)$$

To simplify notation, let us define $T = t+1$ and $M = m+1$. Now, using (31) to approximate $b(T-j, MT, \rho)$ in (29), and approximating sum by an integral, we obtain

$$P_t(0^{m+1}) = \int_1^T \frac{x}{T} \frac{1}{\sqrt{2\pi MT\rho(1-\rho)}} \exp \frac{-(T-x-MT\rho)^2}{2MT\rho(1-\rho)} dx. \tag{32}$$

Integration yields

$$P_t(0^{m+1}) = \sqrt{\frac{M\rho(1-\rho)}{2\pi T}} \left\{ \exp\left(\frac{-(1-T+M\rho T)^2}{2MT\rho(1-\rho)}\right) - \exp\left(\frac{-M\rho T}{2(1-\rho)}\right) \right\} +$$
$$\frac{1}{2}(1-M\rho) \left\{ \mathrm{erf}\left(\frac{M\rho T}{\sqrt{2M\rho(1-\rho)T}}\right) - \mathrm{erf}\left(\frac{1-T+M\rho T}{\sqrt{2M\rho(1-\rho)T}}\right) \right\},$$

where $\mathrm{erf}(x)$ denotes the error function

$$\mathrm{erf}(x) = \frac{2}{\sqrt{\pi}} \int_0^x e^{-t^2} dt. \tag{33}$$

The first term in the above equation (involving two exponentials) tends to 0 with $T \to \infty$. Moreover, since $\lim_{x \to \infty} \mathrm{erf}(x) = 1$, we obtain

$$\lim_{t \to \infty} P_t(0^{m+1}) = \frac{1}{2}(1-M\rho) \left\{ 1 - \lim_{T \to \infty} \mathrm{erf}\left(\frac{1-T+M\rho T}{\sqrt{2M\rho(1-\rho)T}}\right) \right\}.$$

Now, noting that

$$\lim_{T \to \infty} \mathrm{erf}\left(\frac{1-T+M\rho T}{\sqrt{2M\rho(1-\rho)T}}\right) = \begin{cases} 1, & \text{if } M\rho \geq 1, \\ -1, & \text{otherwise,} \end{cases} \tag{34}$$

and returning to the original notation, we recover eq. (25):

$$\lim_{t \to \infty} P_t(0^{m+1}) = \begin{cases} 1-(m+1)\rho & \text{if } \rho < 1/(m+1), \\ 0 & \text{otherwise.} \end{cases} \tag{35}$$